\documentclass[12pt]{amsart}
\usepackage{amsxtra}
\usepackage{amssymb, amsmath}
\usepackage{amsthm}
\usepackage {hyperref}

\newtheorem{lemma}{Lemma}
\newtheorem{theorem}{Theorem}

\newtheorem{proposition}{Proposition}

\title[AG -- Divergence Measures]{GENERALIZED ARITHMETIC AND GEOMETRIC
MEAN DIVERGENCE MEASURE AND THEIR STATISTICAL ASPECTS}

\author{Inder Jeet Taneja}
\address{Inder Jeet Taneja \\
Departamento de Matem\'{a}tica\\ Universidade Federal
de Santa Catarina\\
88.040-900 Florian\'{o}polis, SC, Brazil}

\email{taneja@mtm.ufsc.br}

\urladdr{http://www.mtm.ufsc.br/$\sim $taneja}



\begin{document}

\begin{abstract}
Using Blackwell's definition of comparing two experiments, a
comparison is made with \textit{generalized AG - divergence}
measure having one and two scalar parameters. Connection of
\textit{generalized AG - divergence} measure with \textit{Fisher
measure of information} is also presented. A unified
\textit{generalization of AG - divergence }and\textit{
Jensen-Shannon divergence measures} is also presented.
\end{abstract}

\maketitle

\section{Introduction}

Several measures have been introduced in the literature on
information theory and statistics as measures of information. The
most famous in the literature of statistics is Fisher \cite{fis}
\textit{measure of information}. It measures the amount of
information supplied by data about an unknown parameter $\theta $.
The most commonly used in information theory is the Shannon
\cite{sha} \textit{entropy}. It gives the amount of uncertainty
concerning the outcome of an experiment. Kullback and Leibler
\cite{kul} introduced a measure associated with two distributions
of an experiment. It expresses the amount of information supplied
by the data for discriminating among the distribution. As a
symmetric measure, Jeffreys-Kullback-Leibler \textit{J-divergence}
is commonly used. R\'{e}nyi \cite{ren} generalized both Shannon
\textit{entropy} and Kullback-Leibler \textit{relative
information} by introducing a scalar parameter. Burbea and Rao
\cite{bur1}, \cite{bur2} and Taneja \cite{tan2}, \cite{tan3} have
proposed various alternative ways to generalize the
Jeffreys-Kullback-Leibler \textit{J-divergence}. The proposed
measures of Burbe and Rao \cite{bur1}, \cite{bur2} involve one
parameter. Measures proposed by author \cite{tan3} involve two
scalar parameters.

Let ${\mathcal E} _X = \left\{ {X,S_X ,P_\theta ;\theta \in \Theta
} \right\}$ denote a statistical experiment in which a random
variable or random vector $X$ defined on some sample space $S_X $
is to be observed and the distribution $P_\theta $ of $X$ depends
on the parameter $\theta $ whose values are unknown and lie in
some parameter space $\Theta $. We shall assume that there exists
a generalized probability density function $p(x\vert \theta )$ for
the distribution $P_\theta $ with respect to $\sigma - $finite
measure $\mu $. Let also $\Xi $ denote the class of all prior
distributions on the parameter space $\Theta $. Given a prior
distribution $\xi \in \Xi $, let $p(x)$ denote the corresponding
marginal generalized probability density function (gpdf)
\[
p(x) = \int_\Theta {p(x\vert \theta )d\xi (\theta )} .
\]

Similarly, if we have two prior distributions $\xi _1 $, $\xi _2
\,  \in \Xi $, the corresponding marginal gpdf's are
\[
p_i (x) = \int_\Theta {p(x\vert \theta )d\xi _i (\theta )} , \, i
= 1,2.
\]

In this context, the \textit{relative information}, the
\textit{J--divergence}, the \textit{Jensen-Shannon divergence},
and the \textit{arithmetic and geometric mean divergence}
measures are given as follows:\\

\textbf{$\bullet$ Relative information} (Kullback and Leibler
\cite{kul})
\begin{equation}
\label{eq1} K(\xi _1 ;\xi _2 \vert \vert X) = \int {p_1 (x)\ln
\left( {\frac{p_1 (x)}{p_2 (x)}} \right)d\mu } .
\end{equation}

\bigskip
\textbf{$\bullet$ J -- divergence }(Jeffreys \cite{jef}, Kullback
and Leibler \cite{kul})
\begin{equation}
\label{eq2} J(\xi _1 ;\xi _2 \vert \vert X) = \int {p_1 (x)\ln
\left( {\frac{p_1 (x)}{p_2 (x)}} \right)d\mu } + \int {p_2 (x)\ln
\left( {\frac{p_2 (x)}{p_1 (x)}} \right)d\mu } .
\end{equation}

\bigskip
\textbf{$\bullet$ Jensen--Shannon divergence }(Sibson \cite{sib},
Burbea and Rao \cite{bur1}, \cite{bur2})
\begin{align}
I(\xi _1 ;\xi _2 \vert \vert X) = \frac{1}{2}& \int {\left[ {p_1
(x)\ln p_1 (x) + p_2 (x)\ln p_2 (x)} \right]}\label{eq3}\\
 & - \int {\left( {\frac{p_1 (x) + p_2 (x)}{2}} \right)\ln \left( {\frac{p_1
(x) + p_2 (x)}{2}} \right)d\mu }.\notag
\end{align}

\bigskip
\textbf{$\bullet$ AG -- divergence }(Taneja \cite{tan3})
\begin{equation}
\label{eq4} T(\xi _1 ;\xi _2 \vert \vert X) = \int {\left(
{\frac{p_1 (x) + p_2 (x)}{2}} \right)\ln \left( {\frac{p_1 (x) +
p_2 (x)}{2\sqrt {p_1 (x)p_2 (x)} }} \right)d\mu }.
\end{equation}

\bigskip
The three divergence measures given above can be written in terms
of Kullback-Leibler \textit{relative information} as
\begin{align}
\label{eq5} J(\xi _1 ;\xi _2 \vert \vert X) & = K(\xi _1 ;\xi _2
\vert \vert X) + K(\xi _2 ;\xi _1 \vert \vert X),\\
\label{eq6} I(\xi _1 ;\xi _2 \vert \vert X) & = \frac{1}{2}\left[
{K(\xi _1 ;\frac{\xi _1 + \xi _2 }{2}\vert \vert X) + K(\xi _2
;\frac{\xi _1 + \xi _2 }{2}\vert \vert X)} \right]\\
\intertext{and} \label{eq7} T(\xi _1 ;\xi _2 \vert \vert X) & =
\frac{1}{2}\left[ {K(\frac{\xi _1 + \xi _2 }{2};\xi _1 \vert \vert
X) + K(\frac{\xi _1 + \xi _2 }{2};\xi _2 \vert \vert X)} \right].
\end{align}

Moreover we have the following equality holding among the three
divergence measures
\begin{equation}
\label{eq8} I(\xi _1 ;\xi _2 \vert \vert X) + T(\xi _2 ;\xi _1
\vert \vert X) = \frac{1}{4}J(\xi _1 ;\xi _2 \vert \vert X).
\end{equation}

Recently, author \cite{tan5} proved an interesting inequality
among these three divergence measures:
\begin{equation}
\label{eq9} I(\xi _1 ;\xi _2 \vert \vert X) \leqslant
\frac{1}{8}J(\xi _2 ;\xi _1 \vert \vert X) \leqslant T(\xi _1 ;\xi
_2 \vert \vert X),
\end{equation}

\noindent where all the probability distributions involved are
positive.

In view of (\ref{eq8}) the inequalities (\ref{eq9}) can be
extended as follows:
\begin{equation}
\label{eq10} I(\xi _1 ;\xi _2 \vert \vert X) \leqslant
\frac{1}{8}J(\xi _2 ;\xi _1 \vert \vert X) \leqslant T(\xi _1 ;\xi
_2 \vert \vert X) \leqslant \frac{1}{4}J(\xi _2 ;\xi _1 \vert
\vert X).
\end{equation}

Based on above notations, the Csisz\'{a}r \cite{csi} $\phi -
$\textit{divergence} is given by
\begin{equation}
\label{eq11} C_\phi (\xi _1 ;\xi _2 \vert \vert X) =
\int_{\mathcal X} {p_2 (x)\phi } \left( {\frac{p_1 (x)}{p_2 (x)}}
\right)d\mu ,
\end{equation}

\noindent where the function $\phi $ is arbitrary convex function
defined in the interval $(0,\infty )$. In order to avoid
meaningless expressions, the functions $\phi $ satisfy some
conventional conditions given in \cite{csi}.

In this paper we shall present two parameter generalizations of
the \textit{AG -- divergence}. Also we shall present one
parametric unified generalization of the measures (\ref{eq3}) and
(\ref{eq4}). For two parametric generalization of the measures
(\ref{eq2}) and (\ref{eq3}) refer to Taneja \cite{tan3}. Also
refer on line book by author \cite{tan4}. Here, in this paper we
shall make connections of \textit{generalized AG -- divergence
measures} with Fisher measure of information. The comparison of
experiments is also studied applying Blackwell's \cite{bla}
approach.

\section{Unified $(r,s) - $Arithmetic and Geometric Mean Divergence Measures }

In this section, we shall present two different ways of
generalizing the \textit{AG -- divergence } measure (\ref{eq4}).
Before it we shall give two parametric \textit{unified $(r,s) -
$generalization} \cite{tan3} of the \textit{relative information}:
\begin{equation}
\label{eq12} {\mathcal K}_r^s (\xi _1 ;\xi _2 \vert \vert X) =
\begin{cases}
 {K_r^s (\xi _1 ;\xi _2 \vert \vert X),} & {r \ne 1,\mbox{ }s \ne 1} \\
 {K_1^s (\xi _1 ;\xi _2 \vert \vert X),} & {r = 1,\mbox{ }s \ne 1} \\
 {K_r^1 (\xi _1 ;\xi _2 \vert \vert X),} & {r \ne 1,\mbox{ }s = 1} \\
 {K(\xi _1 ;\xi _2 \vert \vert X),} & {r = 1,\mbox{ }s = 1} \\
\end{cases}
\end{equation}

\noindent for all $r > 0$ and $ - \infty < s < \infty $, where
\[
K_r^s (\xi _1 ;\xi _2 \vert \vert X) = (s - 1)^{ - 1}\left[
{\left( {\int {p_1 (x)^rp_2 (x)^{1 - r}d\mu } } \right)^{\frac{s -
1}{r - 1}} - 1} \right], \, r \ne 1, \, s \ne 1
\]
\[
K_1^s (\xi _1 ;\xi _2 \vert \vert X) = (s - 1)^{ - 1}\left[ {e^{(s
- 1)K(\xi _1 ;\xi _2 \vert \vert X)} - 1} \right], \, s \ne 1
\]

\noindent and
\[
K_r^1 (\xi _1 ;\xi _2 \vert \vert X) = (r - 1)^{ - 1}\ln \left(
{\int {p_1 (x)^rp_2 (x)^{1 - r}d\mu } } \right), \, r \ne 1, \, s
\ne 1.
\]

\subsection*{2.1. First Generalizations}

In (\ref{eq7}) replace $K(\xi _1 ;\xi _2 \vert \vert X)$ by
${\mathcal K}_r^s (\xi _1 ;\xi _2 \vert \vert X)$, we get
\begin{align}
^1{\mathcal T}_r^s (\xi _1 ;\xi _2 \vert \vert X) & =
\frac{1}{2}\left[ {{\mathcal K}_r^s (\frac{\xi _1 + \xi _2
}{2};\xi _1 \vert \vert X) + {\mathcal K}_r^s (\frac{\xi _1 + \xi
_2 }{2};\xi _2 \vert \vert X)} \right] \label{eq13}\\
 & = \begin{cases}
 {^1T_r^s (\xi _1 ;\xi _2 \vert \vert X),} & {r \ne 1,\mbox{ }s \ne 1} \\
 {^1T_1^s (\xi _1 ;\xi _2 \vert \vert X),} & {r = 1,\mbox{ }s \ne 1} \\
 {^1T_r^1 (\xi _1 ;\xi _2 \vert \vert X),} & {r \ne 1,\mbox{ }s = 1} \\
 {T(\xi _1 ;\xi _2 \vert \vert X),} & {r = 1,\mbox{ }s = 1} \\
\end{cases}\notag
\end{align}

\noindent where
\begin{align}
^1T_r^s (\xi _1 ;\xi _2 \vert \vert X) = \,\,& (s - 1)^{ -
1}\frac{1}{2}\left[ {\left\{ {\int {\left( {\frac{p_1 (x) + p_2
(x)}{2}} \right)} } \right.^rp_1 (x)^{1 - r}dx} \right]^{\frac{s -
1}{r - 1}}\notag\\
& + \left. {\left[ {\int {\left( {\frac{p_1 (x) + p_2 (x)}{2}}
\right)} ^rp_2 (x)^{1 - r}d\mu } \right]^{\frac{s - 1}{r - 1}} -
2} \right\} , \, r \ne 1, \, s \ne 1\notag
\end{align}
\[
^1T_1^s (\xi _1 ;\xi _2 \vert \vert X) = \frac{1}{2}(s - 1)^{ -
1}\left[ {e^{(s - 1)K(\frac{\xi _1 + \xi _2 }{2};\xi _1 \vert
\vert X)} + e^{(s - 1)K(\frac{\xi _1 + \xi _2 }{2};\xi _2 \vert
\vert X)} - 2} \right], \, s \ne 1
\]

\noindent and
\begin{align}
T_r^1 (\xi _1 ;\xi _2 \vert \vert X) = &\,\,(r - 1)^{ -
1}\frac{1}{2}\left[ {\ln \left( {\int {\left( {\frac{p_1 (x) + p_2
(x)}{2}} \right)^rp_1 (x)^{1 - r}d\mu } } \right)} \right.\notag\\
& + \left. {\ln \left( {\int {\left( {\frac{p_1 (x) + p_2 (x)}{2}}
\right)^rp_2 (x)^{1 - r}d\mu } } \right)} \right], \, r \ne
1\notag
\end{align}

\noindent for all $r > 0$ and $ - \infty < s < \infty $

\subsection*{2.2. Second Generalizations}

In particular, when $r = s$ in (\ref{eq13}), we get
\begin{align}
\label{eq15} ^1T_s^s (\xi _1 ;\xi _2 \vert \vert X) =& \,\, (s -
1)^{ - 1} \cdot\\
& \cdot \left[ {\int {\left( {\frac{p_1 (x) + p_2 (x)}{2}}
\right)} ^s\left( {\frac{p_1 (x)^{1 - s} + p_2 (x)^{1 - s}}{2}}
\right)d\mu - 1} \right], \notag
\end{align}
\noindent for all $s \ne 1, \, s > 0$.

We shall use the expression (\ref{eq15}) to give the alternative
generalizations of \textit{AG -- divergence} measure. This unified
way is given by
\begin{equation}
\label{eq16} ^2{\mathcal T}_r^s (\xi _1 ;\xi _2 \vert \vert X) =
\begin{cases}
 {^2T_r^s (\xi _1 ;\xi _2 \vert \vert X),} & {r \ne 1,\mbox{ }s \ne 1} \\
 {^2T_1^s (\xi _1 ;\xi _2 \vert \vert X),} & {r = 1,\mbox{ }s \ne 1} \\
 {^2T_r^1 (\xi _1 ;\xi _2 \vert \vert X),} & {r \ne 1,\mbox{ }s = 1} \\
 {T(\xi _1 ;\xi _2 \vert \vert X),} & {r = 1,\mbox{ }s = 1} \\
\end{cases}
\end{equation}

\noindent for all $r > 0$ and $ - \infty < s < \infty $, where
\begin{align}
^2T_r^s (\xi _1 ;\xi _2 \vert \vert X) = & \,\,(s - 1)^{ - 1}
\cdot\notag\\
& \cdot \left\{ {\left[ {\int {\left( {\frac{p_1 (x) + p_2
(x)}{2}} \right)} ^r\left( {\frac{p_1 (x)^{1 - r} + p_2 (x)^{1 -
r}}{2}} \right)d\mu } \right]^{\frac{s - 1}{r - 1}} - 1}
\right\},\notag
\end{align}
\[
^2T_1^s (\xi _1 ;\xi _2 \vert \vert X) = (s - 1)^{ - 1}\left[
{e^{(s - 1)T(\xi _1 ;\xi _2 \vert \vert X)} - 1} \right],
\]

\noindent and
\begin{align}
^2T_r^1 (\xi _1 ;\xi _2 \vert \vert X) = & \,\,(r - 1)^{ - 1}
\cdot\notag\\
& \cdot \ln \left\{ {\int {\left( {\frac{p_1 (x) + p_2 (x)}{2}}
\right)} ^r\left( {\frac{p_1 (x)^{1 - r} + p_2 (x)^{1 - r}}{2}}
\right)d\mu } \right\},\notag
\end{align}
\noindent for all $ r > 0, \,r \ne 1, \, s \ne 1$.

    In particular, we have
\[
    ^1{\mathcal T}_s^s (\xi _1 ;\xi _2 \vert \vert X)
    = \, ^2{\mathcal T}_s^s (\xi _1 ;\xi _2 \vert \vert X).
\]

\subsection*{2.3. Composition Relations}

We observe that the measures $^\alpha {\mathcal T}_r^s (\xi _1
;\xi _2 \vert \vert X)$ ($\alpha = 1,2)$ are continuous with
respect to the parameters $r$ and $s$. This allows us to write
them in the following simplified way
\begin{equation}
\label{eq17} ^\alpha {\mathcal T}_r^s (\xi _1 ;\xi _2 \vert \vert
X) = CE\left\{ {^\alpha T_r^s (\xi _1 ;\xi _2 \vert \vert X)\left|
{r > 0,\mbox{ }r \ne 1,\mbox{ }s \ne 1} \right.} \right\}, \,
\alpha = 1,2,
\end{equation}

\noindent where ``\textit{CE}'' stands for ``\textit{continuous
extension}'' with respect to $r$ and $s$.

Also we can write
\begin{equation}
\label{eq18} ^1{\mathcal T}_r^s (\xi _1 ;\xi _2 \vert \vert X) =
{\mathcal N}_s \left( {^1{\mathcal T}_r^1 (\xi _1 ;\xi _2 \vert
\vert X)} \right)
\end{equation}

\noindent and
\begin{equation}
\label{eq19} ^2{\mathcal T}_r^s (\xi _1 ;\xi _2 \vert \vert X) =
\,\, {\mathcal N}_s \left( {{\mathcal K}_r^1 \left( {\frac{\xi _1
+ \xi _2 }{2};\xi _1 \vert \vert X} \right)} \right)
 + {\mathcal N}_s \left( {{\mathcal K}_r^1 \left( {\frac{\xi _1 + \xi _2 }{2};\xi _2
\vert \vert X} \right)} \right),
\end{equation}

\noindent where ${\mathcal N}_s :(0,\infty ) \to
\mathbb{R}$(reals) is given by
\begin{equation}
\label{eq20} {\mathcal N}_s (x) = \begin{cases}
 {(s - 1)^{ - 1}\left[ {e^{(s - 1)x} - 1} \right],} & {s \ne 1} \\
 {x,} & {s = 1} \\
\end{cases}
\end{equation}

\begin{proposition} \label{pro21} The measure ${\mathcal N}_s (x)$ given
above has the following properties:

\begin{itemize}
\item[(i)] ${\mathcal N}_s (x) \geqslant 0$ with equality iff $x =
0$;

\item[(ii)] ${\mathcal N}_s (x)$ is an increasing function of $x$;

\item[(iii)] ${\mathcal N}_s (x)$ is an increasing function of
$s$;

\item[(iv)] ${\mathcal N}_s (x)$ is strictly convex function of
$x$ for $s > 1$;

\item[(v)] ${\mathcal N}_s (x)$ is strictly concave function of
$x$ for $s < 1$.
\end{itemize}
\end{proposition}

\subsection*{2.4. Alternative Generalizations}

We see that the measure (\ref{eq15}) is considered for $s > 0$. It
is required for the non-negativity of the measure. We can rewrite
it in little different way, where we don't require this condition.
This form is given by
\begin{align}
\label{eq21} IT_s (\xi _1 ;\xi _2 \vert \vert X) = & \,\,\left[
{s(s - 1)} \right]^{ - 1} \cdot\\
& \cdot \left[ {\int {\left( {\frac{p_1 (x) + p_2 (x)}{2}}
\right)} ^s\left( {\frac{p_1 (x)^{1 - s} + p_2 (x)^{1 - s}}{2}}
\right)d\mu - 1} \right], \notag
\end{align}
\noindent where $\, s \ne 0,1$

The measure (\ref{eq21}) admits the following limiting cases:
\[
\mathop {\lim }\limits_{s \to 0} IT_s (\xi _1 ;\xi _2 \vert \vert X) = I(\xi
_1 ;\xi _2 \vert \vert X)
\]

\noindent and
\[
\mathop {\lim }\limits_{s \to 1} IT_s (\xi _1 ;\xi _2 \vert \vert X) = T(\xi
_1 ;\xi _2 \vert \vert X),
\]

\noindent where $I(\xi _1 ;\xi _2 \vert \vert X)$ and $T(\xi _1
;\xi _2 \vert \vert X)$ are as given by (\ref{eq3}) and
(\ref{eq4}) respectively.

In view of these limiting cases, we re-write the measure
(\ref{eq21}) in the following unified way
\begin{equation}
\label{eq22} {\mathcal {IT}}_s (\xi _1 ;\xi _2 \vert \vert X) =
\begin{cases}
 {IT_s (\xi _1 ;\xi _2 \vert \vert X),} & {s \ne 0,1} \\
 {I(\xi _1 ;\xi _2 \vert \vert X),} & {s = 0} \\
 {T(\xi _1 ;\xi _2 \vert \vert X),} & {s = 1} \\
\end{cases}
\end{equation}

\section{Relationship with Csisz\'{a}r $\phi - $Divergence}

We can relate the above generalizations of the \textit{AG --
divergence} measure with the well known Csisz\'{a}r \textit{$\phi
- $divergence}. It is given as follows:
\[
^1{\mathcal T}_r^s (\xi _1 ;\xi _2 \vert \vert X) =
\frac{1}{2}\left[ {\eta _s \left( {\phi (\xi _1 ;\xi _2 \vert
\vert X)^{\frac{1}{r - 1}}} \right) + \eta _s \left( {\phi ^\ast
(\xi _1 ;\xi _2 \vert \vert X)^{\frac{1}{r - 1}}} \right)} \right]
\]

\noindent and
\[
^2{\mathcal T}_r^s (\xi _1 ;\xi _2 \vert \vert X) = \eta _s \left(
{\phi ^ - (\xi _1 ;\xi _2 \vert \vert X)^{\frac{1}{r - 1}}}
\right),
\]

\noindent where
\[
\eta _s (y) = \begin{cases}
 {(s - 1)^{ - 1}\left[ {y^{s - 1} - 1} \right],} & {s \ne 1} \\
 {y,} & {s = 1} \\
\end{cases}
\]

\noindent and $\phi (\xi _1 ;\xi _2 \vert \vert X)$, $\phi ^\ast
(\xi _1 ;\xi _2 \vert \vert X)$ and $\phi ^ - (\xi _1 ;\xi _2
\vert \vert X)$ are the $\phi - $\textit{divergences} of $\xi _1
$, $\xi _2 $ in the notations of Vajda \cite{vaj} with
\[
\phi (x) = \left( {\frac{1 + x}{2}} \right)^r,
\]
\[
\phi ^\ast (x) = x\phi \left( {\frac{1}{x}} \right)
\]

\noindent and
\[
\phi ^ - (x) = \frac{1}{2}\left[ {\phi (x) + \phi ^\ast (x)} \right].
\]

We can also write the measure (\ref{eq22}) in terms of Csisz\'{a}r
$\phi - $\textit{divergence} as follows:
\begin{equation}
\label{eq23} {\mathcal {IT}}_s (\xi _1 ;\xi _2 \vert \vert X) =
\int_{\mathcal X} {p_2(x)\phi _{{\mathcal {IT}}_s } \left(
{\frac{p_1(x)}{p_2(x)}} \right)} d\mu ,
\end{equation}

\noindent where
\begin{equation}
\label{eq24} \phi _{{\mathcal {IT}}_s } (x) = \begin{cases}
 {\left[ {2s(s - 1)} \right]^{ - 1}\left[ {\left( {x^{1 - s} + 1}
\right)\left( {\frac{x + 1}{2}} \right)^s - (x + 1)} \right],} & {s \ne 0,1}
\\
 {\frac{x}{2}\ln x + \left( {\frac{x + 1}{2}} \right)\ln \left( {\frac{2}{x
+ 1}} \right),} & {s = 0} \\
 {\left( {\frac{x + 1}{2}} \right)\ln \left( {\frac{x + 1}{2\sqrt x }}
\right),} & {s = 1} \\
\end{cases},
\end{equation}

\noindent for all $x \in (0,\infty )$.

\begin{proposition} \label{pro31} For all $r > 0$ and $ - \infty < s <
\infty $, we have

\begin{itemize}
 \item[(i)] $^\alpha {\mathcal T}_r^s (\xi _1 ;\xi _2
\vert \vert X) \geqslant 0 \,\,(\alpha = 1,2)$;
 \item[(ii)] $
^1{\mathcal T}_r^s (\xi _1 ;\xi _2 \vert \vert X)\begin{cases}
 { \leqslant \mbox{ }^2{\mathcal T}_r^s (\xi _1 ;\xi _2 \vert \vert X),} & {s
\leqslant r} \\
 { \geqslant \mbox{ }^2{\mathcal T}_r^s (\xi _1 ;\xi _2 \vert \vert X),} & {s
\geqslant r} \\
\end{cases}$

\item[(iii)] ${\mathcal {IT}}_s (\xi _1 ;\xi _2 \vert \vert X)
\geqslant 0 $.
\end{itemize}
\end{proposition}

\section{Divergence Measures and Sufficiency of Experiments}

Blackwell \cite{bla} definition of comparison of experiments
states that experiment ${\mathcal E} _X $ is sufficient for
experiment ${\mathcal E} _Y $, denoted by ${\mathcal E} _X \succeq
{\mathcal E} _Y $, if there exists a stochastic transformation of
$X$ to a random variable $Z(X)$ such that for each $\theta \in
\Theta $ the random variable $Z(X)$ and $Y$ have identical
distributions. By ${\mathcal E} _Y = \left\{ {Y,S_Y ,Q_\theta
;\theta \in \Theta } \right\}$ we shall denote a second
statistical experiment for which there exists a \textit{gpdf}
$g(y\vert \theta )$ for the distribution $Q$ with respect to a
$\sigma - $finite measure $\mu $. According to this definition, if
${\mathcal E} _X \succeq {\mathcal E} _Y $, then there exists a
nonnegative function $h$ satisfying (DeGroot \cite{deg})
\begin{equation}
\label{eq25} g(y\vert \theta ) = \int_{\mathcal X} {h(y\vert
x)f(x} \vert \theta )d\mu
\end{equation}

\noindent and
\[
\int_{\mathcal X} {h(y\vert x)} d\upsilon = 1.
\]

Changing the order of integration in (\ref{eq25}), we get
\begin{equation}
\label{eq26} g_i (y ) = \int_{\mathcal X} {h(y\vert x)f_i (x}
)d\mu , \, i = 1,2.
\end{equation}

Let $I$ be any measure of information contained in an experiment.
If ${\mathcal E} _X \succeq {\mathcal E} _Y $ implies that $I_X
\geqslant I_Y $, then we say that ${\mathcal E} _X $ is as
informative as ${\mathcal E} _Y $. This approach is successfully
carried out by Lindley \cite{lin} for Shannon \textit{entropy}.
Goel and DeGroot \cite{god} applied it for Kullback and Leibler
\cite{kul} \textit{relative information}. Ferentinos and
Papaioannou \cite{fep} applied for $\alpha - $\textit{order
generalization} of Kullback and Leibler \textit{relative
information} and generalizations of Fisher \textit{measure of
information}. Author \cite{tan1} extended it to different
generalizations of \textit{J--divergence} measure having two
scalar parameters. For the \textit{I -- divergence measure} and
their two parametric generalizations refer to Taneja et al.
\cite{tpm}. Here our aim is to compare experiments for the
\textit{unified }$(r,s) - $\textit{AG -- divergences }given by
(\ref{eq13}) and (\ref{eq16}). Results are also extended for the
measure (\ref{eq22}).

\begin{theorem} \label{the41} If ${\mathcal E} _X \succeq {\mathcal E} _Y $, then
$^\alpha {\mathcal T}_r^s (\xi _1 ;\xi _2 \vert \vert X) \geqslant
\, ^\alpha {\mathcal T}_r^s (\xi _1 ;\xi _2 \vert \vert Y)$
($\alpha = 1,2)$ for every $\xi _1 $, $\xi _2  \, \in \Xi $, for
all $r > 0$ and $ - \infty < s < \infty $.
\end{theorem}

\begin{proof} Since ${\mathcal E} _X \succeq {\mathcal E} _Y $, there exists a
function $h$ satisfying (\ref{eq25}) and (\ref{eq26}), then we can
write
\begin{align}
\label{eq27} & \left( {\frac{g_1 (y) + g_2 (y)}{2}} \right)^r g_1
(y)^{1 - r}\\
 & = \left[ {\int_{\mathcal X} {h(y\vert x)\left( \frac{f_1 (x)+ f_2 (x)} {2} \right)d\mu } } \right]^{r}\left[ {\int_{\mathcal X} {h(y\vert
x)f_1 (x)d\mu } } \right]^{1-r}.\notag
\end{align}

Applying H\"{o}lder's inequality on the right side of
(\ref{eq27}), we get
\begin{align}
\label{eq28} & \left( {\frac{g_1 (y) + g_2 (y)}{2}} \right)^rg_1
(y)^{1 - r}\\
& \begin{cases}
 { \geqslant \int_{\mathcal X} {\left[ {h(y\vert x)\left( {\frac{f_1 (x) + f_2 (x)}{2}}
\right)} \right]^{r}} \left[ {h(y\vert x)f_1 (x)}
\right]^{1-r}d\mu ,} & {0 < r <
1} \\\\
 { \leqslant \int_{\mathcal X} {\left[ {h(y\vert x)\left( {\frac{f_1 (x) + f_2 (x)}{2}}
\right)} \right]^{r}} \left[ {h(y\vert x)f_1 (x)}
\right]^{1-r}d\mu ,} & {r > 1}
\\
\end{cases}\notag
\end{align}

Hence
\begin{align}
\label{eq29} \int_{\mathcal Y} & {\left( {g_1 (y)}
\right)^{1-r}\left(
{\frac{g_1 (y) + g_2 (y)}{2}} \right)^{r}d\upsilon }\\
& \begin{cases}
 { \geqslant \int_{\mathcal X} {\left( {f_1 (x)} \right)^{1-r}} \left( {\frac{f_1 (x) + f_2
(x)}{2}} \right)^{r}d\mu ,} & {0 < r < 1} \\\\
 { \leqslant \int_{\mathcal X} {\left( {f_1 (x)} \right)^{1-r}} \left( {\frac{f_1 (x) + f_2
(x)}{2}} \right)^{r}d\mu ,} & {r > 1} \\
\end{cases}\notag
\end{align}

As $sign\left( {\frac{s - 1}{r - 1}} \right) = sign(s - 1)$ for $r
> 1$ and $sign\left( {\frac{s - 1}{r - 1}} \right) \ne sign(s -
1)$ for $0 < r < 1$, where $sign(x) = 1$ if $x > 0$ and $sign(x) =
- 1$ if $x < 0$, then from (\ref{eq29}) one gets
\begin{align}
\label{eq30} \frac{1}{s - 1}& \left[ {\int_{\mathcal Y} {\left(
{g_1 (y)} \right)^{1-r}\left( {\frac{g_1 (y) + g_2 (y)}{2}}
\right)^{r}d\upsilon } } \right]^{\frac{s - 1}{r - 1}}\\
& \leqslant \frac{1}{s - 1}\left[ {\int_{\mathcal X} {\left( {f_1
(x)} \right)^{1-r}\left( {\frac{f_1 (x) + f_2 (x)}{2}}
\right)^{r}d\mu } } \right]^{\frac{s - 1}{r - 1}},\notag
\end{align}

\noindent for all $r \ne 1$, $s \ne 1$, $r > 0$.

Similarly, we can obtain
\begin{align}
\label{eq31} \frac{1}{s - 1}& \left[ {\int_{\mathcal Y} {\left(
{g_2 (y)} \right)^{1-r}\left( {\frac{g_1 (y) + g_2 (y)}{2}}
\right)^{r}d\upsilon } } \right]^{\frac{s - 1}{r - 1}}\\
 & \leqslant \frac{1}{s - 1}\left[ {\int_{\mathcal X} {\left( {f_2 (x)} \right)^{1-r}\left(
{\frac{f_1 (x) + f_2 (x)}{2}} \right)^{r}d\mu } } \right]^{\frac{s
- 1}{r - 1}},\notag
\end{align}

\noindent for all $r \ne 1$, $s \ne 1$, $r > 0$

Adding (\ref{eq30}) and (\ref{eq31}), subtracting $2(s - 1)^{ -
1}$ ($s \ne 1)$ and then dividing by 2, we get
\[
^1T_r^s (\xi _1 ;\xi _2 \vert \vert X) \geqslant \,^1T_r^s (\xi _1
;\xi _2 \vert \vert Y), \, r \ne 1, \, s \ne 1, \, r
> 0.
\]

Since the unified measure $^1{\mathcal T}_r^s (\xi _1 ;\xi _2
\vert \vert X)$ given in (\ref{eq13}) is a continuous extension of
$^1T_r^s (\xi _1 ;\xi _2 \vert \vert X)$ for the real parameters
$r$ and $s$ we can immediately conclude that
\[
^1{\mathcal T}_r^s (\xi _1 ;\xi _2 \vert \vert X) \geqslant
\,^1{\mathcal T}_r^s (\xi _1 ;\xi _2 \vert \vert Y), \, r > 0,
\]

\noindent whenever ${\mathcal E} _X \succeq {\mathcal E} _Y $.\\

Let us prove now the second part. Since ${\mathcal E} _X \succeq
{\mathcal E} _Y $, there exist a function $h$ satisfying
(\ref{eq25}) and (\ref{eq26}), then we can write
\begin{align}
\label{eq32} & \left( {\frac{g_1 (y) + g_2 (y)}{2}}
\right)^r\left( {\frac{g_1 (y)^{1 - r} + g_2 (y)^{1 - r}}{2}}
\right)\\
 & =  \frac{1}{2}\left[ {\int_{\mathcal X} {h(y\vert x)f_1 (x)d\mu } } \right]^{1-r}\left[
{\int_{\mathcal X} {h(y\vert x)\left( {\frac{f_1 (x) + f_2
(x)}{2}} \right)d\mu } } \right]^{r}\notag\\
& \,\,\,\,+ \frac{1}{2}\left[ {\int_{\mathcal X} {h(y\vert x)f_2
(x)d\mu } } \right]^{1-r}\left[ {\int_{\mathcal X} {h(y\vert
x)\left( {\frac{f_1 (x) + f_2 (x)}{2}} \right)d\mu } }
\right]^{r}.\notag
\end{align}

Applying H\"{o}lder's inequality in (\ref{eq32}), integrating over
$\mathcal Y$, and using the fact that $\int_{\mathcal Y} {h(y\vert
x)d\upsilon } = 1$, we get
\begin{align}
\label{eq33} \int_{\mathcal Y} & {\left( {\frac{g_1 (y) + g_2
(y)}{2}} \right)^r\left( {\frac{g_1 (y)^{1 - r} + g_2 (y)^{1 -
r}}{2}} \right)d\upsilon }\\
& \begin{cases}
 { \geqslant \int_{\mathcal X} {\left( {\frac{f_1 (x) + f_2 (x)}{2}} \right)^r\left(
{\frac{f_1 (x)^{1 - r} + f_2 (x)^{1 - r}}{2}} \right)d\mu ,} } & {0 < r < 1}
\\\\
 { \leqslant \int_{\mathcal X} {\left( {\frac{f_1 (x) + f_2 (x)}{2}} \right)^r\left(
{\frac{f_1 (x)^{1 - r} + f_2 (x)^{1 - r}}{2}} \right)d\mu ,} } & {r > 1} \\
\end{cases}.\notag
\end{align}

As $sign\left( {\frac{s - 1}{r - 1}} \right) = sign(s - 1)$ for $r
> 1$ and $sign\left( {\frac{s - 1}{r - 1}} \right) \ne sign(s -
1)$ for $0 < r < 1$, we have
\begin{align}
\label{eq34} \frac{1}{s - 1}& \left[ {\int_{\mathcal Y} {\left(
{\frac{g_1 (y) + g_2 (y)}{2}} \right)^r\left( {\frac{g_1 (y)^{1 -
r} + g_2 (y)^{1 - r}}{2}} \right)d\upsilon } } \right]^{\frac{s -
1}{r - 1}}\\
& \leqslant \frac{1}{s - 1}\left[ {\int_{\mathcal X} {\left(
{\frac{f_1 (x) + f_2 (x)}{2}} \right)^r\left( {\frac{f_1 (x)^{1 -
r} + f_2 (x)^{1 - r}}{2}} \right)d\mu } } \right]^{\frac{s - 1}{r
- 1}}.\notag
\end{align}

Subtracting $(s - 1)^{ - 1}$ ($s \ne 1)$ on both sides of
(\ref{eq34}), we get
\[
^2T(\xi _1 ,\xi _2 \vert \vert X) \geqslant \, ^2T(\xi _1 ,\xi _2
\vert \vert Y), \, r \ne 1, \, s \ne 1, \, r > 0,
\]

\noindent and consequently, we have
\[
^2\mathcal T_r^s (\xi _1 ,\xi _2 \vert \vert X) \geqslant \,
^2\mathcal T_r^s (\xi _1 ,\xi _2 \vert \vert Y), \, r > 0,
\]

\noindent whenever ${\mathcal E} _X \succeq {\mathcal E} _Y $.
\end{proof}

\begin{theorem} \label{the42} If ${\mathcal E} _X \succeq {\mathcal E} _Y $,
then ${\mathcal {IT}}_s (\xi _1 ;\xi _2 \vert \vert X) \geqslant
{\mathcal {IT}}_s (\xi _1 ;\xi _2 \vert \vert Y)$ for every $\xi
_1 $, $\xi _2  \, \in \Xi $, for all $s \in ( - \infty ,\infty )$.
\end{theorem}

The proof of the above theorem is based on the following lemmas.

\begin{lemma} \label{lem41} (Joint convexity). If $\phi :(0,\infty )
\to \mathbb{R}$ be convex, then $C_\phi (\xi _1 ;\xi _2 \vert
\vert X)$ jointly convex for every $\xi _1 $, $\xi _2  \,  \in \Xi
$. Moreover if $\phi (1) = 0$, then $C_\phi (\xi _1 ;\xi _2 \vert
\vert X) \geqslant 0$.
\end{lemma}

\begin{lemma} \label{lem42} If ${\mathcal E} _X \succeq {\mathcal E} _Y $, then
$C_\phi (\xi _1 ;\xi _2 \vert \vert X) \geqslant C_\phi (\xi _1
;\xi _2 \vert \vert Y)$ for every $\xi _1 $, $\xi _2  \, \in \Xi
$, provided $\phi $ is convex.
\end{lemma}

Lemma \ref{lem41} is due to Csisz\'{a}r \cite{csi} and Lemma
\ref{lem42} is due to Ferentinos and Papaioannou \cite{fep}.

\begin{proof} \textit{of Theorem \ref{the41}}. In view of Lemmas \ref{lem41}
and \ref{lem42}, it is sufficient to prove the convexity of the
function $\phi _{{\mathcal {IT}}_s } (x)$ given by (\ref{eq19}).
It is in view of the following derivatives:
\begin{equation}
\label{eq35} \phi _{{\mathcal {IT}}_s } ^\prime (x) =
\begin{cases}
 {\frac{s\left( {\frac{x + 1}{2x}} \right)^{1 - s} + (1 - s)\left(
{\frac{x^s + 1}{2}} \right)\left( {\frac{x + 1}{2}} \right)^{ - s}}{2s(s -
1)},} & {s \ne 0,1} \\
 {\frac{1}{4}\left[ {1 - x^{ - 1} - \ln x - 2\ln \left( {\frac{2}{x + 1}}
\right)} \right],} & {s = 0} \\
 {\frac{1}{2}\left[ {\ln x + \ln \left( {\frac{2}{x + 1}} \right)} \right],}
& {s = 1} \\
\end{cases}
\end{equation}

\noindent and
\begin{equation}
\label{eq36} \phi _{{\mathcal {IT}}_s } ^{\prime \prime }(x) =
\begin{cases}
 {\frac{1}{8}\left( {x^{s - 2} + 1} \right)\left( {\frac{x + 1}{2}}
\right)^{ - s - 1},} & {s \ne 0,1} \\
 {\frac{1}{4}\left( {\frac{x^2 + 1}{x^2(x + 1)}} \right),} & {s = 0} \\
 {\frac{1}{2x(x + 1)},} & {s = 1} \\
\end{cases}.
\end{equation}

Thus we have $\phi _{{\mathcal {IT}}_s } ^{\prime \prime }(x) > 0$
for all $x > 0$, and hence, $\phi _{{\mathcal {IT}}_s } (x)$ is
\textit{convex} for all $x > 0$. Also, we have $\phi _{{\mathcal
{IT}}_s } (1) = 0$. In view of this we can say that \textit{I--
$\&$ T-- divergence of type s} given by (\ref{eq22}) is
\textit{nonnegative} and \textit{convex} in the pair of
probability distributions $P$ and $Q$.
\end{proof}

\section{$\phi - $Divergence and Fisher Information Matrix}

Consider a family ${\mathcal M} = \left\{ {{\mathcal P}_\theta
,\theta \in \Theta } \right\}$ of probability measures on a
measurable space $(X,{\mathcal A})$ dominated by a finite or
$\sigma - $finite measure $\mu $. The parameter space $\Theta $
can either be an open subset of the real line or an open subset of
$n - $dimensional Euclidean space $\mathbb{R}^k$. Let $f(x,\theta
) = \frac{dP_\theta }{d\mu }$. Let $\Gamma = \left\{ {f(x,\theta
)\left| {x \in X,\mbox{ }\theta \in \Theta } \right.} \right\}$.

The Fisher \cite{fis} measure of information is given by
\begin{equation}
\label{eq37} I_X^F (\theta ) = \begin{cases}
 {E_\theta \left[ {\frac{\partial }{\partial \theta }\ln f(x,\theta )}
\right]^2,} & {\mbox{if }\theta \mbox{ is univariate}} \\
 {E_\theta \left\| {\frac{\partial }{\partial \theta _i }\ln f(x,\theta
)\frac{\partial }{\partial \theta _j }\ln f(x,\theta )} \right\|_{k\times k}
,} & {\mbox{if }\theta \mbox{ is }k - \mbox{variate}} \\
\end{cases}
\end{equation}

\noindent where $\left\|(\cdot) \right\|_{k\times k} $ denotes a
$k\times k$ matrix and $E_\theta $ denotes the expectation with
respect to $f(x,\theta )$, where $f(x,\theta ) \in \Gamma $. Let
us suppose that the following regularity conditions are
satisfied:\\

\begin{itemize}

\item[(a)] $\frac{\partial }{\partial \theta _i }f(x,\theta )$
exists for all $x \in X$, all $\theta \in \Theta $, and all $i =
1,2,...,k$.\\

\item[(b)] For any $A \in {\mathcal A}$,

$\frac{d}{d\theta _i }\int_A {f(x,\theta )d\mu = \int_A
{\frac{\partial }{\partial \theta _i }f(x,\theta )d\mu } } ,\,$for
all $i = 1,2,...,k.$\\
\end{itemize}

For $f_1 ,f_2 \in \Gamma $, the Csisz\'{a}r \cite{csi} $\phi -
$\textit{divergence} can be re-written as
\begin{equation}
\label{eq38} K_\phi (f_1 \vert \vert f_2 ) = \int {f_2 \phi \left(
{\frac{f_1 }{f_2 }} \right)d\mu } ,
\end{equation}

\noindent with $\phi (1)$ not necessarily zero and $\phi (x)$ is a
continuously differentiable nonnegative real function. As usual,
the function $\phi (x)$ is generally supposed to be convex, but
here we don't assume that $\phi (x)$ is convex.

Following Kagan \cite{kag} and Ferentinos and Papaioannou
\cite{fep}, we define
\begin{equation}
\label{eq39}
I_{ij}^C (\theta ) = \mathop {\lim }\limits_{t \to 0} \inf
\frac{1}{t^2}\left\{ {K_\phi \left( {f(x,\theta )\vert \vert
\frac{f(x,\theta + te_i ) + f(x,\theta + te_j )}{2}} \right) - \phi (1)}
\right\}.
\end{equation}

Then the Csisz\'{a}r parametric matrix is given by
\begin{equation}
\label{eq40} I_X^C (\theta ) = \left\| {I_{ij}^C (\theta )}
\right\|_{k\times k} ,
\end{equation}

\noindent where $\theta + te_i ,\theta + te_j \in \Theta $, $i,j =
1,2,...,k$ and $e_1 (1,0,...,0)$, $e_2 (0,1,...,0)$, ..., $e_k
(0,0,...,1)$ are the unit vectors.\\

Suppose the following conditions hold:\\
\begin{itemize}
\item[(c)] $\int {\left| {\frac{\partial ^2}{\partial \theta _i
\partial \theta _j }f(x,\theta )d\mu } \right|} < \infty $ for all
$\theta \in \Theta $ and $i,j = 1,2,...,k$.\\

\item[(d)] The third order partial derivative of $f(x,\theta )$
with respect to $\theta $ exists for all $\theta \in \Theta $ and
$x \in X.$
\end{itemize}

Based on the above considerations the following theorem holds.

\begin{theorem} \label{the51} If the conditions (a)-(d) are satisfied,
then for all $\theta \in \Theta $, we have
\begin{equation}
\label{eq41} I_X^C (\theta ) = \begin{cases}
 {\frac{{\phi }''(1)}{2}I_X^F (\theta ),} & {\mbox{if }\theta \mbox{ is
univariate}} \\\\
 {\frac{{\phi }''(1)}{2}\left[ {S_X^F (\theta ) + I_X^F (\theta )} \right],}
& {\mbox{if }\theta \mbox{ is }k - \mbox{variate}} \\
\end{cases}
\end{equation}

\noindent where
\[
S_X^F (\theta ) = \frac{1}{2}\left[ {M_X^F (\theta ) + M_X^F (\theta )^T}
\right]
\]

\noindent with
\[
M_X^F (\theta ) = \left\| {{\begin{array}{*{20}c}
 {I_{11}^F (\theta )} \hfill & {I_{11}^F (\theta )} \hfill & \cdots \hfill &
{I_{11}^F (\theta )} \hfill \\
 {I_{22}^F (\theta )} \hfill & {I_{22}^F (\theta )} \hfill & \cdots \hfill &
{I_{22}^F (\theta )} \hfill \\
 \vdots \hfill & \vdots \hfill & \ddots \hfill & \vdots \hfill \\
 {I_{kk}^F (\theta )} \hfill & {I_{kk}^F (\theta )} \hfill & \cdots \hfill &
{I_{kk}^F (\theta )} \hfill \\
\end{array} }} \right\| = \left\| {I_{ij}^F (\theta )} \right\|_{k\times k}
\]

\noindent and
\[
I_{ij}^F (\theta ) = E_\theta \left[ {\frac{\partial }{\partial \theta _i
}\ln f(x,\theta )} \right]^2.
\]
\end{theorem}

This result has been derived by Aggarwal \cite{agg}. Similar
results derived for the R\'{e}nyi, Kagan, Kullback-Leibler,
Matusita measures of information can be seen in Kagan \cite{kag},
Aggarwal \cite{agg}, Boekee \cite{boe}, Ferentinos and Papaioannou
\cite{fep}, Taneja \cite{tan1}, Salicr\'{u} and Taneja \cite{sat},
etc.

\section{Unified $(r,s)$--T--Divergence and Fisher
Information Matrix}

To get the relationship between \textit{unified
$(r,s)$-T-divergence} and Fisher information matrix, first we give
the following proposition due to Salicr\'{u} and Taneja
\cite{sat}.

\begin{proposition} \label{pro51} Let
\[
G_\phi ^h (f_1 \vert \vert f_2 ) = h\left( {K_\phi (f_1 \vert \vert f_2 ) -
\phi (1)} \right),
\]

\noindent where $h$ is a continuous differentiable real function
with $h(0) = 0$, and $K_\phi (f_1 \vert \vert f_2 )$ is given by
(\ref{eq38}). Suppose the conditions (a)-(d) are satisfied. Then
for $\theta \in \Theta $, we have
\[
G_\phi ^h \left( {I_X^C (\theta )} \right) = {h}'(0)I_X^C (\theta ),
\]

\noindent where
\[
G_\phi ^h \left( {I_X^C (\theta )} \right) = \left\| {G_\phi ^h \left(
{I_{ij}^C (\theta )} \right)} \right\|_{k\times k}
\]

\noindent with

\begin{align}
G_\phi ^h & \left( {I_X^C (\theta )} \right) \notag\\
& = \mathop {\lim }\limits_{t \to 0} \inf \frac{1}{t^2}h\left(
{K_\phi \left( {f(x,\theta )\vert \vert \frac{f(x,\theta + te_i )
+ f(x,\theta + te_j )}{2}} \right) - \phi (1)} \right)\notag
\end{align}

\noindent and $I_X^C (\theta )$ is the Csisz\'{a}r information
matrix given in (\ref{eq40}).
\end{proposition}

Now we shall apply the above results to connect the measures
(\ref{eq13}), (\ref{eq16}) and (\ref{eq22}) with Fisher measure of
information.

\begin{proposition} \label{pro52} If the conditions (a)-(d)
are satisfied, then for all $\theta \in \Theta $, we have
\begin{align}
^1{\mathcal T}_r^s (\theta ) & = \frac{r}{8}\left[ {S_X^F (\theta
) + I_X^F (\theta )} \right],\\
^2{\mathcal T}_r^s (\theta ) & = \frac{r}{8}\left[ {S_X^F (\theta
) + I_X^F (\theta )} \right]\\
\intertext{and} {\mathcal {IT}}_s (\theta ) & = \frac{1}{8}\left[
{S_X^F (\theta ) + I_X^F (\theta )} \right].
\end{align}
\end{proposition}

\begin{proof} We shall prove for each part separately.\\

We can write
\[
^1T_r^s (f_1 \vert \vert f_2 ) = h\left( {K_{\phi _1 } (f_1 \vert \vert f_2
) - \phi _1 (1)} \right) + h\left( {K_{\phi _2 } (f_1 \vert \vert f_2 ) -
\phi _2 (1)} \right),
\]

\noindent where
\[
\phi _1 (x) = x\left( {\frac{1 + x}{2x}} \right)^r, \, r \ne 1, \,
r > 0
\]
\[
\phi _2 (x) = \left( {\frac{1 + x}{2}} \right)^r, \, r \ne 1, \, r
> 0
\]

\noindent and
\[
h(x) = \left[ {2(s - 1)} \right]^{ - 1}\left[ {(x + 1)^{\frac{s -
1}{r - 1}} - 1} \right], \, r \ne 1, \, s \ne 1, \, r > 0.
\]

This gives
\[
\phi _1 ^{\prime \prime }(1) = \phi _2 ^{\prime \prime }(1) = \frac{r(r -
1)}{4}
\]

\noindent and
\[
{h}'(0) = \left[ {2(r - 1)} \right]^{ - 1}, \, r \ne 1, \, r > 0
\]

We have
\[
^1T_r^s (\theta ) = \frac{r}{16}\left[ {S_X^F (\theta ) + I_X^F
(\theta )} \right], \, r \ne 1, \, r > 0
\]

\noindent and consequently,
\[
^1{\mathcal T}_r^s (\theta ) = \frac{r}{16}\left[ {S_X^F (\theta )
+ I_X^F (\theta )} \right],
\]

\noindent for all $\theta \in \Theta $, $0 < r < \infty $ and $ -
\infty < r < \infty $.\\

Again, we can write
\[
^2T_r^s (f_1 \vert \vert f_2 ) = h\left( {K_\phi (f_1 \vert \vert f_2 ) -
\phi _1 (1)} \right),
\]

\noindent where
\[
\phi (x) = \left( {\frac{1 + x}{2}} \right)^r\left( {\frac{x^{1 -
r} + 1}{2}} \right), \, r \ne 1, \, r > 0
\]

\noindent and
\[
h(x) = (s - 1)^{ - 1}\left[ {(x + 1)^{\frac{s - 1}{r - 1}} - 1}
\right], \, r \ne 1, \, s \ne 1, \, r > 0.
\]

This gives
\[
{\phi }''(1) = \frac{r(r - 1)}{4}
\]

\noindent and
\[
{h}'(0) = (r - 1)^{ - 1}, \, r \ne 1, \, r > 0
\]

We have
\[
^2T_r^s (\theta ) = \frac{r}{16}\left[ {S_X^F (\theta ) + I_X^F
(\theta )} \right], \, r \ne 1, \, r > 0
\]

\noindent and consequently,
\[
^2{\mathcal T}_r^s (\theta ) = \frac{r}{16}\left[ {S_X^F (\theta )
+ I_X^F (\theta )} \right],
\]

\noindent for all $\theta \in \Theta $, $0 < r < \infty $ and $ -
\infty < s < \infty $.\\

It is easy to check that $\phi _{{\mathcal {IT}}_s } ^{\prime
\prime }(1) = \frac{1}{4}$ for $s \in ( - \infty ,\infty )$. This
gives
\[
{\mathcal {IT}}_s (\theta ) = \frac{1}{8}\left[ {S_X^F (\theta ) +
I_X^F (\theta )} \right],
\]

\noindent for all $\theta \in \Theta $.
\end{proof}

\end{document}